\documentclass[a4paper,10pt]{amsart}
\usepackage[arrow,matrix]{xy}
\usepackage{amsmath,amssymb,amscd,bbm,amsthm,mathrsfs}
\usepackage[colorlinks=true,citecolor=blue,linkcolor=blue]{hyperref}

\theoremstyle{plain}
\newtheorem{thm}{Theorem}[section]
\newtheorem{cor}[thm]{Corollary}
\newtheorem{lem}[thm]{Lemma}
\newtheorem{prop}[thm]{Proposition}

\newtheorem{rem}[thm]{Remark}

   \def\op{\oplus} \def\ot{\otimes}
\def\Hom{\operatorname {Hom}}

\def\Ext{\operatorname {Ext}} 
 \def\k{\mathbbm{k}}

\begin{document}
\title{\bf Deformations of Koszul Artin-Schelter Gorenstein algebras}

\author{Ji-Wei He, Fred Van Oystaeyen and Yinhuo Zhang}
\address{J.-W. He\newline \indent Department of Mathematics, Shaoxing College of Arts and Sciences, Shaoxing Zhejiang 312000,
China\newline \indent Department of Mathematics and Computer
Science, University of Antwerp, Middelheimlaan 1, B-2020 Antwerp,
Belgium} \email{jwhe@usx.edu.cn}
\address{F. Van Oystaeyen\newline\indent Department of Mathematics and Computer
Science, University of Antwerp, Middelheimlaan 1, B-2020 Antwerp,
Belgium} \email{fred.vanoystaeyen@ua.ac.be}
\address{Y. Zhang\newline
\indent Department WNI, University of Hasselt, Universitaire Campus,
3590 Diepenbeek, Belgium} \email{yinhuo.zhang@uhasselt.be}

\date{}
\begin{abstract}
We compute the Nakayama automorphism of a PBW-deformation of a Koszul Artin-Schelter Gorenstein algebra of finite global dimension, and give a criterion for an augmented PBW-deformation of a Koszul Calabi-Yau algebra to be Calabi-Yau. The relations between the Calabi-Yau property of augmented PBW-deformations and that of non-augmented cases are discussed. The Nakayama automorphisms of PBW-deformations of Koszul Artin-Schelter Gorenstein algebras of global dimensions 2 and 3 are given explicitly. We show that if a PBW-deformation of a graded Calabi-Yau algebra is still Calabi-Yau, then it is defined by a potential under some mild conditions. Some classical results are also recovered. Our main method used in this paper is elementary and based on linear algebra. The results obtained in this paper will be applied in a subsequent paper \cite{HVZ1}.
\end{abstract}


\maketitle

\section*{Introduction}

Dubois-Violette showed in \cite{DV0,DV} that a (generalized) Koszul Artin-Schelter (AS, for short) Gorenstein algebra of finite global dimension is determined by a multi-linear form, and hence, to some extent, the study of the AS-Gorenstein property of a Koszul algebra is equivalent to dealing with certain problems in linear algebra. Motivated by this observation, we ask if we can convert homological problems of Koszul AS-Gorenstein algebras to linear algebra problems. The first question arising is whether it is possible to write down explicitly the {\it Nakayama automorphism} (for the definition, see the following section) of an AS-Gorenstein algebra in a linear algebra way.

Instead of using Dubois-Violette's multi-linear form, we write the generating relations of a Koszul algebra by a sequence of matrices. It turns out that this is an efficient way to discuss certain homological problems. For example, by using these matrices, we may write down the product of the Yoneda Ext-algebra of the Koszul algebra explicitly; then we write down the Nakayama automorphism of a Koszul AS-Gorenstein algebra by matrices through Van den Bergh's formula \cite{VdB} (and later, Berger and Marconnet \cite{BM}) for the Nakayama automorphism of a Koszul AS-Gorenstein algebra. We also recover the result that graded Calabi-Yau algebras of dimension 3 are obtained from a superpotential by some simple calculations with matrices, compared to the result in \cite{Bo} where Bocklandt proved the result by a sophisticated analysis of the minimal projective resolution of the trivial module.

The aim of this paper is to study the Poincar\'{e}-Birkhoff-Witt (PBW) deformations of AS-Gorentein algebras. We try to write down the Nakayama automorphisms of PBW-deformations of Koszul AS-Gorenstein algebras, in particular of AS-Gorenstein algebras of global dimension 2 or 3, and then study the Calabi-Yau property of PBW-deformations.

The paper is organized as follows. In Section 1, we set up some terminology and notation used in this paper, and recall some definitions and basic properties of $N$-Koszul algebras and AS-Gorenstein algebras. In Section 2, we discuss the Nakayama automorphisms of PBW-deformations of general Koszul AS-Gorenstein algebras of finite global dimensions. The main results of this section are Proposition \ref{mprop2}, Theorems \ref{mthm2} and \ref{mthm3}. Wu and Zhu in \cite{WZ} proved a result similar to Theorem \ref{mthm2} under the assumption that the associated graded algebra is Noetherian. In contrast, we do not need the Noetherian condition. Also our method is totally different from that of \cite{WZ}. Theorem \ref{mthm3} extends a result we proved in \cite{HVZ}.

Section 3 is devoted to compute the Nakayama automorphisms of AS-Gorenstein algebras of global dimension 2. The main result of this section is Theorem \ref{mthm1}. Some results obtained in \cite{Ber} and \cite{DV} are recovered by some linear computations.

Section 4 is devoted to discuss the PBW-deformations of AS-Gorenstein algebras of global dimension 3. It is well known that an AS-Gorenstein algebra of global dimension 3 is $N$-Koszul. We write the generating relations of a Koszul algebra $A$ by a sequence of matrices, allowing to write down the product of the homogeneous dual algebra $A^!$ of $A$ explicitly. Similarly, we also write down the Yoneda product of $E(A)=\op_{i=0}^3\Ext_A^i(\k_A,\k_A)$ when $A$ is an $N$-Koszul AS-Gorenstein algebra of global dimension 3. The Nakayama automorphisms of PBW-deformations are expressed in Proposition \ref{mprop13} by some matrices. Some simple computations on the matrices yield that a connected graded Calabi-Yau algebra of dimension 3 is defined by a superpotential, hence we recover Bocklandt's result \cite{Bo}. Moreover, we prove that if a PBW-deformation of a graded Calabi-Yau algebra of dimension 3 is still Calabi-Yau, then it is defined by a potential in case that the PBW-deformation is augmented or the associated graded algebra is a domain (Theorems \ref{mthm5}). Moreover, the potential for a Calabi-Yau PBW-deformation can be explicitly written down.

\section{Preliminaries}

Throughout $\k$ is a field of characteristic zero. All the vector spaces, algebras and coalgebras involved in this paper are over $\k$. We write $\ot$ for $\ot_\k$. Let $M=\op_{i\in \mathbb{Z}}M_i$ be a graded vector space. For an integer $n$, we write $M(n)$ for the graded vector space whose $i^{th}$ component is $M(n)_i=M_{i+n}$. Let $B$ be an algebra, $\phi$ be an automorphism of $B$. We write ${}_1B_\phi$ for the $B$-bimodule whose left $B$-action is  the regular action and right $B$-action is twisted by $\phi$. We write $B^e$ for the enveloping algebra $B\ot B^{op}$.

Let $V$ be a finite dimensional vector space, and let $T(V)=\k\op V\op V^{\ot 2}\op\cdots$ be the tensor algebra. Given a subspace $R\subseteq V^{\ot N}$ ($N\ge2$), the graded algebra $A=T(V)/(R)$ is called an $N$-homogeneous algebra, where $(R)$ is the ideal of $T(V)$ generated by $R$. Let $V^*=\Hom_\k(V,\k)$, and $R^\perp\subseteq (V^*)^{\ot N}$ be the orthogonal complement of $R$. The $N$-homogeneous algebra $A^!=T(V^*)/(R^\perp)$ is called the {\it homogeneous dual} of $A$. Since $A$ is graded, we may view the ground field as a graded left (or right) $A$-module. For any finitely generated graded left (or right) $A$-module $M$, the minimal graded projective resolution of $M$ exists. Suppose that the minimal graded projective resolution of $\k_A$ is the following: $$\cdots\longrightarrow P^{-i}\longrightarrow\cdots\longrightarrow P^{-1}\longrightarrow P^0\longrightarrow\k_A\longrightarrow0.$$ If, for all $i\ge0$, the graded module $P^{-i}$ is generated in degree $\kappa(i)$, then $A$ is called an {\it $N$-Koszul algebra} \cite{Ber1}, where $\kappa:\mathbb{N}\to \mathbb{N}$ is a function defined by $\kappa(i)=\frac{i}{2}N$ in case $i$ is even, or $\kappa(i)=\frac{i-1}{2}N+1$ in case $i$ is odd. In particular, if $N=2$, then $A$ is called a {\it Koszul algebra} \cite{Pr}. For a Koszul algebra $A$, its homogeneous dual algebra $A^!$ is also a Koszul algebra. Moreover, the Yoneda algebra $E(A)=\op_{i\ge0}\Ext_A^i(\k_A,\k_A)$ is isomorphic to $A^!$ \cite{BGS}. In general, the homogeneous dual algebra of an $N$-Koszul algebra $A$ is no longer an $N$-Koszul algebra if $N\ge3$, and in this case, $\Ext^i_A(\k_A,\k_A)\cong A^!_{\kappa(i)}$. We refer to \cite{Pr,BGS,PPo} for further properties of Koszul algebras, and \cite{Ber1,GMMZ} for $N$-Koszul algebras if $N\ge3$.

Let $A=A_0\op A_1\op\cdots$ be a graded algebra. A {\it PBW-deformation} of $A$ is a filtered algebra $U$ with an ascending filtration $0\subseteq F_0U\subseteq F_1U\subseteq F_2U\subseteq\cdots$ such that the associated graded algebra $gr(U)$ is isomorphic to $A$.

{\bf Convention:} henceforth, by a {\it deformation} we mean a PBW-deformation of a graded algebra.

If $A=T(V)/(R)$ is a 2-homogeneous algebra, then a deformation $U$ of $A$ is determined by two linear maps $\nu:R\to V$ and $\theta:R\to\k$ in sense that $U\cong T(V)/(r-\nu(r)-\theta(r):r\in R)$. The linear maps $\nu$ and $\theta$ satisfy some Jacobian type conditions (more details, see\cite{BG,PPo}). If $\theta=0$, then $U$ is called an {\it augmented} deformation of $A$.

Generally, if $A=T(V)/(R)$ is an $N$-homogeneous algebra, then a deformation $U$ of $A$ is determined by $N$ linear maps $\alpha_i:R\to V^{\ot N-i}$ for $i=1,\dots,N$, where we set $V^{\ot 0}=\k$, in sense that $U\cong T(V)/(r+\alpha_1(r)+\cdots+\alpha_N(r):r\in R)$. Also, the linear maps $\alpha_i$'s satisfy certain Jacobian type conditions (see \cite{BGi} or \cite{FV}). If $\alpha_N=0$, then $U$ is called an augmented deformation of $A$.

We are interested in the deformations of $N$-Koszul AS-Gorenstein algebras. Let $A=\k\op A_1\op A_2\op\cdots$ be a connected graded algebra. Assume that $A$ has global dimension $d<\infty$. $A$ is called an {\it AS-Gorenstein algebra} if $\Ext_A^i(\k_A,A)=0$ if $i\neq d$ and $\Ext_A^d(\k_A,A)\cong \k$. Compared to the classical definition, we do not assume that an AS-Gorenstein algebra is Noetherian.

An AS-Gorenstein algebra can be viewed as a kind of generalization of a Frobenius algebra. Let $B=\k\op B_1\op B_2\op\cdots$ be a finite dimensional graded algebra. Recall that $B$ is called a {\it graded Frobenius algebra} if there is an isomorphism of graded left $B$-modules: $$\Theta:B\cong B^*(-d).$$ In this case, we have $B_d\neq 0$ and $B_i=0$ for all $i>d$. We say that the {\it length} of the Frobenius algebra $B$ is $d$. In general, the isomorphism $\Theta$ is not a $B$-bimodule morphism. In fact, there is a unique graded automorphism $\phi$ of $B$ such that $\Theta:{}_1B_\phi\to B^*$ is an isomorphism of $B$-bimodules. The isomorphism $\phi$ is called the {\it Nakayama automorphism} of $B$. If $\phi=id$, then $B$ is called a {\it symmetric} algebra. The isomorphism $\Theta$ induces a nondegenerate bilinear form $\langle\ ,\ \rangle:B\times B\to \k$ defined by $\langle a,b\rangle=\Theta(1)(ab)$ for all $a,b\in B$. We call a Frobenius algebra $B$ of length $d$ {\it graded symmetric} if $\langle a,b\rangle=(-1)^{i(d-i)}\langle b,a\rangle$ for all $a\in B_i$ and $b\in B_{d-i}$.

We have the following relation between $N$-Koszul AS-Gorenstein algebras and graded Frobenius algebras.

\begin{lem}\label{lem11} Let $A$ be an $N$-Koszul algebra of global dimension $d$. Then $A$ is AS-Gorenstein if and only if $E(A)=\op_{i\ge0}\Ext^i_A(\k_A,\k_A)$ is a graded Frobenius algebra.
\end{lem}
When $A$ is Koszul, the proof of the lemma above can be found in \cite{Smi}; when $A$ is a general $N$-Koszul algebra, the proof can be found in \cite{BM}, more generally, see \cite{LPWZ}.

\begin{lem} \label{lem12} Let $A=T(V)/(R)$ be an $N$-Koszul AS-Gorenstein algebra of global dimension $d$. Then $\Ext^i_{A^e}(A,A\ot A)=0$ for $i\neq d$, and $\Ext_{A^e}^d(A,A\ot A)\cong {}_1A_\zeta(\kappa(d))$ for some graded automorphism $\zeta$ of $A$.
\end{lem}

The graded automorphism $\zeta$ is unique. We call $\zeta$ the {\it Nakayama automorphism} of $A$.

The proof of the lemma above can be found in \cite{VdB} when $A$ is Koszul and Noetherian, and in the proof of \cite[Theorem 6.3]{BM} when $A$ is a general $N$-Koszul algebra. Moreover, the Nakayama automorphism $\zeta$ of $A$ is also constructed in \cite{VdB} and \cite{BM}. In fact, by Lemma \ref{lem11}, the Yoneda algebra $E(A)$ is graded Frobenius. Assume that $\phi$ is the Nakayama automorphism of $E(A)$. Since $\Ext_A^1(\k_A,\k_A)\cong V^*$, the automorphism $\phi$ induces a graded automorphism $\varphi$ of $A$. Let $\varepsilon$ be the automorphism of $A$ defined by multiplying $(-1)^n$ on a homogeneous element $a\in A_n$. Then the Nakayama automorphism $\zeta$ in Lemma \ref{lem12} is $$\zeta=\varepsilon^{d+1}\phi^{-1}.$$

We end this section by recalling the definition of a Calabi-Yau algebra. An algebra $B$ is called a {\it Calabi-Yau} algebra of dimension $d$  \cite{Gin}, if (i) $B$ is
homologically smooth, that is; $B$ has a bounded resolution of
finitely generated projective $B$-bimodules, (ii)
$\Ext^i_{B^e}(B,B\ot B)=0$ if $i\neq d$ and $\Ext_{B^e}^d(B,B\ot B)\cong
B$ as $B$-bimodules. Note that if, further, $B$ is a connected graded algebra, then $B$ must be an AS-Gorenstein algebra \cite{BT}.

\section{Deformations of Koszul AS-Gorenstein algebras}\label{gendef}

Let $A=T(V)/(R)$ be a Koszul algebra, and let $U=T(V)/(r-\nu(r)-\theta(r):r\in R)$ be a deformation of $A$. Assume $\dim V=n$ and $\{x_1,\dots,x_n\}$ be a basis of $V$. Let $\{x_1^*,\dots,x_n^*\}$ be the dual basis of $V^*$. Define $$C_{-m}=(R\ot V^{\ot m-2})\cap (V\ot R\ot V^{\ot m-3})\cap\cdots\cap (V^{\ot m-2}\ot R)$$ for $m\ge2$, $C_0=\k$ and $C_{-1}=V$. Let $C=\op_{m\ge0} C_{-m}$. Then $C$ is a graded subcoalgebra of the tensor coalgebra $T(V)$. The map $\nu:R\to V$ induces a graded coderivation $\delta_C$ of degree 1 on the coalgebra $C$. We extend the linear map $\theta:R\to \k$ to a graded map (also denoted by $\theta$) $\theta:C\to \k$ of degree 2 in an obvious way. Since $U$ is a deformation of $A$, we see that the graded coalgebra $C$ together with the coderivation $\delta_C$ and the graded map $\theta$ form a {\it curved differential graded coalgebra} in sense of Positselski \cite{Po}. Let $A^!=T(V^*)/(R^\perp)$ be the dual algebra of $A$. The dual map $\nu^*:V^*\to R^*$ induces a derivation $\delta_{A^!}$ of degree 1 on $A^!$. If we view the map $\theta:R\to \k$ as an element in $A^!_2$, then the graded algebra $A^!$ together with the derivation $\delta_{A^!}$ and the element $\theta$ form a {\it curved differential graded algebra}, that is, $\delta_{A^!}$ is a derivation of $A^!$ and $\theta\in A^!_2$ such that $\delta_{A^!}^2(\alpha)=\theta\alpha-\alpha\theta$ \cite{PPo,Po}. The triple $(A^!,\delta_{A^!},\theta)$ is sometimes called the curved differential graded algebra {\it dual} to $U$. The graded dual algebra of the curved differential graded coalgebra $(C,\delta_C,\theta)$ is isomorphic to the curved differential graded algebra $(A^!,\delta_{A^!},\theta)$ \cite[Section 6]{Po}. So, we can view $C$ as a graded $A^!$-bimodule.

We construct a sequence of $U$-bimodules by using the curved differential graded coalgebra $(C,\delta_C,\theta)$ as follows: {\small\begin{equation}\label{mtag6}
    \cdots\longrightarrow U\ot C_{-m}\ot U\overset{D}\longrightarrow\cdots\longrightarrow U\ot C_{-2}\ot U\overset{D}\longrightarrow U\ot C_{-1}\ot  U\overset{D}\longrightarrow U\ot U\overset{\mu}\longrightarrow U\longrightarrow0,
\end{equation}}
where the last morphism $\mu$ is the multiplication of $U$, and the morphism $D$ is defined as, for $c\in C_{-m}$, $a,b\in U$ and $m\ge1$,
$$D(a\ot c\ot b)=\sum_{i=1}^nax_i\ot c\cdot x_i^*\ot b+(-1)^m\sum_{i=1}^na\ot x_i^*\cdot c\ot x_ib-a\ot \delta_C(c)\ot b.$$
Using the equations in \cite[Proposition 1.1, Ch. 5]{PPo}, we can check that $D^2=0$. The complex (\ref{mtag6}) coincides with the complex constructed by the natural twisting cochain from $C$ to $U$ as shown in \cite[Section 6]{Po}. We endow each component of the complex (\ref{mtag6}) with an ascending filtration by setting $$\begin{array}{l}
               F_i(U\ot C_{-m}\ot U)=0\text{ for $i<m$},\\
               F_m(U\ot C_{-m}\ot U)=F_0U\ot C_{-m}\ot F_0U,\text{ and}\\
               \displaystyle F_{m+k}(U\ot C_{-m}\ot U)=\sum_{i+j=k}F_iU\ot C_{-m}\ot F_jU, \text{ for $k\ge1$}.
             \end{array}
$$
Then the differential $D$ of the complex (\ref{mtag6}) is compatible with the filtration. The first level of the spectral sequence induced by the filtration is exactly the classical two-sided Koszul resolution of $A$. Hence the complex (\ref{mtag6}) is exact, and then it is a projective resolution of the bimodule ${}_UU_U$.

Applying the functor $\Hom_{U^e}(-,U\ot U)$ to the projective resolution of the bimodule ${}_UU_U$ and observing that $A^!$ is dual to $C$, we get the following complex:
{\small\begin{equation}\label{mtag4}
    0\longrightarrow U\ot A^!_0\ot U\overset{\widehat{D}}\longrightarrow U\ot A^!_1\ot U\overset{\widehat{D}}\longrightarrow\cdots\longrightarrow U\ot A^!_{m-1}\ot U\overset{\widehat{D}}\longrightarrow U\ot A^!_m\ot U\longrightarrow\cdots,
\end{equation}} where the differential $\widehat{D}$ is defined as follows: for any $a,b\in U$ and $\alpha\in A^!_m$,
$$\widehat{D}(a\ot \alpha\ot b)=\sum_{i=1}^na\ot x^*_i\alpha\ot x_ib+(-1)^{m+1}\sum_{i=1}^nax_i\ot \alpha x_i^*\ot b-a\ot \delta_{A^!}(\alpha)\ot b.$$
The complex (\ref{mtag4}) inherits the filtration of the complex (\ref{mtag6}) as follows: $$F_i(U\ot A^!_m\ot U)=0,\text{ for $i<-m$},$$ and $$F_{-m+k}(U\ot A^!_m\ot U)=\sum_{i+j=k}F_iU\ot A_m^!\ot F_jU, \text{ for $k\ge0$}.$$
It is easy to see that the differential $\widehat{D}$ preserves the filtration.

Now assume that $A=T(V)/(R)$ is a Koszul AS-Gorenstein algebra of global dimension $d$, and that $U=T(V)/(r-\nu(r)-\theta(r):r\in R)$ is a deformation of $A$. The projective resolution of the bimodule ${}_UU_U$ constructed above reads as follows: $$0\longrightarrow U\ot C_{-d}\ot U\overset{D}\longrightarrow\cdots\longrightarrow U\ot C_{-2}\ot U\overset{D}\longrightarrow U\ot C_{-1}\ot  U\overset{D}\longrightarrow U\ot U\overset{\mu}\longrightarrow U\longrightarrow0.$$ After applying the functor $\Hom_{U^e}(-,U\ot U)$ to this resolution, we obtain
\begin{equation}\label{mtag7}
    0\longrightarrow U\ot A^!_0\ot U\overset{\widehat{D}}\longrightarrow U\ot A^!_1\ot U\overset{\widehat{D}}\longrightarrow\cdots\longrightarrow U\ot A^!_{d-1}\ot U\overset{\widehat{D}}\longrightarrow U\ot A^!_d\ot U.
\end{equation}

We deduce the following result from the complex above. Yekutieli proved a more general result under the assumption that the AS-Gorenstein algebra is Noetherian \cite{Yek}. While, in our case, we do not need the Noetherian condition on $A$.

\begin{lem} \label{lem21} Let $A$ be a Koszul AS-Gorenstein algebra of global dimension $d$, and let $U$ be a deformation of $A$. Assume that $\zeta$ is the Nakayama automorphism of $A$. Then we have $\Ext_{U^e}^i(U,U\ot U)=0$ for $i\neq d$ and $\Ext_{U^e}^d(U,U\ot U)\cong {}_1U_\xi$ as $U$-bimodules, where $\xi$ is a filtration-preserving automorphism of $U$ such that $gr(\xi)=\zeta$.
\end{lem}
\proof Let $E^{pq}_\bullet$ be the spectral sequence obtained from the filtration of the complex (\ref{mtag7}). The first level $E_1^{pq}$ of the spectral sequence is exactly the complex obtained from the Koszul bimodule complex of $A=gr(U)$ (see \cite[Theorem 9.1]{VdB}), and hence collapses by Lemma \ref{lem12}. Since the filtration is ascending and bounded below, the spectral sequence converges. Then we have $\Ext_{U^e}^i(U,U\ot U)=0$ for $i\neq d$ and $\Ext_{U^e}^d(U,U\ot U)\cong W$, where $W$ is a filtered $U$-bimodule such that $gr(W)\cong {}_1A_\zeta$. Now applying \cite[Lemmas 4 and 5]{VdB1}, we obtain the desired result. \qed

As in the graded case, the automorphism $\xi$ in the above lemma is unique up to inner automorphisms. We call $\xi$ the {\it Nakayama automorphism} of $U$.

Let $A$ be a Koszul AS-Gorenstein algebra of global dimension $d$. By Lemma \ref{lem11}, $A^!$ is a graded Frobenius algebra. Hence there is a graded automorphism $\phi:A^!\to A^!$ such that there is an isomorphism of graded $A^!$-bimodules $$\Theta:{}_1A^!_\phi\to {A^!}^*(-d).$$ Assume that $\phi$ is defined by $\phi(x_1^*,\dots,x_n^*)=(x_1^*,\dots,x^*_n)P$, where $P=(p_{ij})$ is an invertible $n\times n$ matrix. The automorphism $\phi$ induces an automorphism $\varphi:A\to A$, which is determined by its restriction to $A_1=V$, namely, $\varphi(x_1,\dots,x_n)=(x_1,\dots,x_n)P^t$.

By Lemma \ref{lem21}, the complex (\ref{mtag7}) is exact except at the final position. The homology at the final position is ${}_1U_\xi$. So, there is a bimodule epimorphism $\pi:U\ot A^!_d\ot U\to {}_1U_\xi$. Let $\varpi\in A^!_d$ be the element such that $\Theta(1)(\varpi)=1$. Since the differential of the complex (\ref{mtag7}) preserves the filtration and the cohomology at the final position of the associated graded complex is ${}_1A_{\varepsilon^{d+1}\varphi^{-1}}$ (Lemma \ref{lem12}), we have $\pi(\varpi)=u\in F_0U=\k$. Since $A^!$ is Frobenius, $\dim A^!_{d-1}=n$. Let $w_1,\dots,w_n\in A^!_{d-1}$ be the elements such that $x_i^*w_j=\delta_j^i\varpi$ for all $i$ and $j$. Then $\{w_1,\dots,w_n\}$ forms a basis of $A^!_{d-1}$, and $w_jx_i^*=\phi^{-1}(x_i^*)w_j$. Assume $\delta_{A^!}(w_j)=\lambda_j\varpi$, $j=1,\dots,n$. For $a,b\in U$, we have $$\pi\circ \widehat{D}(a\ot w_j\ot b)=0.$$ Recall that $\phi$ is defined by the matrix $P$. Thus $\phi^{-1}$ is defined by $P^{-1}$. Let $P^{-1}=(l_{ij})$. Then we have:
$$\begin{array}{ccl}
    0 & = &\pi\circ \widehat{D}(a\ot w_j\ot b)  \\
     & = & \pi(\sum_{i=1}^na\ot x^*_iw_j\ot x_ib+(-1)^{d}\sum_{i=1}^nax_i\ot w_j x_i^*\ot b-a\ot \delta_C(w_j)\ot b) \\
    &  = & \pi(a\ot \varpi\ot x_jb+(-1)^{d}\sum_{i=1}^nax_i\ot \phi^{-1}(x_i^*)w_j\ot b-a\ot \lambda_j\varpi\ot b))  \\
    &=& ua\xi(x_jb)+(-1)^d\sum_{i=1}^nuax_il_{ji}\xi(b)-ua\lambda_j\xi(b).
  \end{array}
$$
In the equations above, if we set $b=1$ and cancel out $u$, then we obtain
\begin{equation}\label{mtag5}
    a\xi(x_j)+(-1)^d\sum_{i=1}^nax_il_{ji}-a\lambda_j=0.
\end{equation}
Since $\xi$ is a filtered automorphism and $gr(\xi)=\varepsilon^{d+1}\varphi^{-1}$ (Lemmas \ref{lem12} and \ref{lem21}), we have $$\xi(x_j)=(-1)^{d+1}\sum_{i=1}^{d+1}l_{ji}x_i+k_j,$$ for some $k_j\in\k$. Now the equation (\ref{mtag5}) is equivalent to $$(-1)^{d+1}\sum_{i=1}^nal_{ji}x_i+ak_i+(-1)^d\sum_{i=1}^nax_il_{ji}-a\lambda_j=0,$$ which, in turn, is equivalent to $$ak_j-a\lambda_j=0,\text{ for all $a\in U$ and $j=1,\dots,n$}.$$ Hence, we get $$k_j=\lambda_j,\text{ for all $j=1,\dots,n$}.$$

Summarizing we obtain the following result.

\begin{prop} \label{mprop2} Let $A=T(V)/(R)$ be a Koszul AS-Gorenstein algebra of global dimension $d$, and let $A^!$ be its dual algebra. Assume that $\{x_1,\dots,x_n\}$ is a basis of $V$, and $\{x_1^*,\dots,x_n^*\}$ is the dual basis of $V^*$.

Let $U=T(V)/(r-\nu(r)-\theta(r):r\in R)$ be a deformation of $A$, and let $(A^!,\delta_{A^!},\theta)$ be the curved differential graded algebra dual to $U$. Assume that the Nakayama automorphism $\phi:A^!\to A^!$ is defined by $\phi(x_1^*,\dots,x_n^*)=(x_1^*,\dots,x_n^*)P$, where $P$ is an $n\times n$ matrix. Choose a basis $\varpi$ of $A^!_d$, and assume that $\{w_1,\dots,w_n\}$ is the basis of $A^!_{d-1}$ such that $x_i^*w_j=\delta^i_j\varpi$. Assume further  $\delta_{A^!}(w_i)=\lambda_i\varpi$ for all $i=1,\dots,n$, and let $\mathbf{\lambda}=(\lambda_1,\dots,\lambda_n)\in \k^n$. Then $\Ext_{U^e}^i(U,U\ot U)=0$ for $i\neq d$, and  $$\Ext^d_{U^e}(U,U\ot U)\cong {}_1U_{\xi},$$ where $\xi:U\to U$ is an automorphism defined by the linear map
\begin{equation}\label{mtag8}
    \chi:V\to V\op \k,\quad\chi(x_1,\dots,x_n)=(-1)^{d+1}(x_1,\dots,x_n){(P^{-1}})^t+\mathbf{\lambda}.
\end{equation}
\end{prop}

If the Koszul algebra $A$ in the proposition is Calabi-Yau, then $A^!$ is a graded symmetric algebra \cite[Proposition A.5.2]{Bo}. Hence in this case the Nakayama automorphism of $A^!$ is $\phi=\epsilon^{d+1}$, where $\epsilon$ is the automorphism of $A^!$ defined as: $\epsilon(f)=(-1)^if$ for $f\in A^!_i$. Then the proposition implies the following result.

\begin{cor} \label{cor1} Let $A$ be a Koszul Calabi-Yau algebra of dimension $d$, and let $U=T(V)/(r-\nu(r)-\theta(r):r\in R)$ be a PBW-deformation of $A$, and let $(A^!,\delta_{A^!},\theta)$ be the curved differential graded algebra dual to $U$. If $\delta_{A^!}(A^!_{d-1})=0$, then $U$ is a Calabi-Yau algebra.
\end{cor}

Let $(A^!,\delta_{A^!},\theta)$ be the curved differential graded algebra dual to $U=T(V)/(r-\nu(r)-\theta(r):r\in R)$. If $\theta$ lies in the center of $A^!$, then $(A^!,\delta_{A^!})$ is actually a differential graded algebra. In this case, $U'=T(V)/(r-\nu(r):r\in R)$ is an augmented deformation of $A$ by \cite[Proposition 4.1]{PPo}. The proposition above implies the following result.

\begin{prop}\label{mprop3}
Let $A$, $U$ and $(A^!,\delta_{A^!},\theta)$ be as in Proposition \ref{mprop2}. If $\theta$ is in the center of $A^!$, then $U'=T(V)/(r-\nu(r):r\in R)$ is an augmented deformation of $A$. Moreover, we have $\Ext^i_{{U'}^e}(U',U'\ot U')=0$ for $i\neq d$, and $$\Ext^d_{{U'}^e}(U',U'\ot U')\cong {}_1U'_{\xi'},$$ where the automorphism $\xi':U'\to U'$ is defined by the same linear map (\ref{mtag8}) as in Proposition \ref{mprop2}.
\end{prop}

In particular, if $A$ is Calabi-Yau, then we have the next result.

\begin{thm}\label{mthm2}
Let $A=T(V)/(R)$ be a Koszul Calabi-Yau algebra of global dimension $d$. Assume that $U$ is an augmented deformation of $A$, and that $(A^!,\delta_{A^!})$ is the differential graded algebra dual to $U$. Then the following are equivalent:

{\rm(i)} $U$ is a Calabi-Yau algebra;

{\rm(ii)} $E(U)=\op_{i=1}^d\Ext^i_U(\k,\k)$ is a graded symmetric algebra of length $d$;

{\rm(iii)} $\delta_{A^!}(A^!_{d-1})=0$.
\end{thm}
\proof (i) $\Longrightarrow$ (ii) is follows from \cite[Proposition A.5.2]{Bo}.

(ii) $\Longrightarrow$ (iii). By the Frobenius property of $E(U)$, we have $\Ext^d_U(\k,\k)\neq 0$. Then by \cite[Proposition 6.1]{PPo}, the $d^{th}$ cohomology of the differential graded algebra $(A^!,\delta_{A^!})$ is not zero. Since $A^!_d$ is of dimension one, we have $\delta_{A^!}(A^!_{d-1})=0$.

(iii) $\Longrightarrow$ (ii). Since $A$ is Calabi-Yau, $A^!$ is graded symmetric of length $d$ by \cite[Proposition A.5.2]{Bo}. Hence the matrix $P$ in Proposition \ref{mprop2} is $P=(-1)^{d+1}I_n$, where $I_n$ is the $n\times n$ unit matrix. That (iii) implies (i) (and therefore (ii)) is Corollary \ref{cor1}. \qed

\begin{rem} {\rm Wu and Zhu proved in \cite{WZ} the equivalence of (i) and (iii) under the assumption that $A$ is Noetherian. We do not need the Noetherian condition. Moreover, the method used in \cite{WZ} is totally different from ours.}
\end{rem}

The theorem above implies the following result, which extends \cite[Theorem 5.3]{HVZ}.

\begin{thm}\label{mthm3} Let $A=T(V)/(R)$ be a Koszul Calabi-Yau algebra of global dimension $d$. Suppose that both $U=T(V)/(r-\nu(r)-\theta(r):r\in R)$ and $U'=T(V)/(r-\nu(r):r\in R)$ are deformations of $A$. If $U'$ is Calabi-Yau, then so is $U$.

Conversely, if $U$ is Calabi-Yau and $A$ is a domain, then $U'$ is Calabi-Yau.
\end{thm}
\proof If $U'$ is Calabi-Yau, then by Theorem \ref{mthm2}(iii), we have $\delta_{A^!}(A^!_{d-1})=0$. Now Proposition \ref{mprop2} implies that the Nakayama automorphism of $U$ is the identity, and hence, $U$ is Calabi-Yau.

Conversely, if $U$ is Calabi-Yau, then the automorphism $\xi$ in Proposition \ref{mprop2} is inner. If $A$ is domain, then the inner automorphism group of $U$ is trivial. Hence $\xi$ is identity map. Therefore $\delta_{A^!}(A^!_{d-1})=0$. Then Proposition \ref{mprop3} implies that $U'$ is Calabi-Yau. \qed

In the following sections, we will apply the results obtained in this section to the deformations of AS-Gorenstein algebras of global dimensions 2 and 3. The Nakayama automorphisms will be computed in detail.

\section{Frobenius algebras of Length 2 and AS-Gorenstein algebras of global dimension 2}

Let $E=\k\op E_1\op E_2$ such that $\dim E_1=n$ and $\dim E_2=1$. The set of graded Frobenius structures on $E$ corresponds to $GL(\k,n)$ in the following way.

Fix a basis of $E_1$, say, $\{x_1,\dots,x_n\}$ and a basis of $E_2$, $\{z\}$. For an invertible $n\times n$-matrix $M\in GL(\k,n)$, we define a multiplication on $E$ as follows: for $x,y\in E_2$, assuming $x=a_1x_1+\cdots a_nx_n$ and $y=b_1x_1+\cdots+b_nx_n$, then $xy=\mathbf{a}M\mathbf{b}^tz$, where $\mathbf{a}=(a_1,\dots,a_n)\in\k^n$ and $\mathbf{b}=(b_1,\dots,b_n)\in \k^n$. Then $E$ is a graded Frobenius algebra with a nondegenerate bilinear form $\langle\ ,\ \rangle$ defined by $\langle x,y\rangle=\mathbf{a}M\mathbf{b}^t$ for all $x,y\in E_1$. A straightforward check shows that $\langle x,y\rangle=\langle y,\phi(x)\rangle$ where $\phi:E\to E$ is a graded algebra automorphism defined by $$\phi(x)=(x_1,\dots,x_n)M^{-1}M^t\mathbf{a}^t\text{ and }\phi(z)=z.$$ Hence $\phi$ is the Nakayama automorphism of $E$.

Conversely, it is easy to see that any graded Frobenius algebra of length 2 is defined in this way.

For further discussion, we need more notation. Let $W$ and $W'$ be $n$-dimensional vector spaces with fixed bases $\{w_1,\dots,w_n\}$ and $\{w'_1,\dots,w'_n\}$ respectively. For an element $\alpha\in W\ot W'$, we may write $\alpha=\sum_{i,j=1}^np_{ij}w_i\ot w'_j$. Let $P=(p_{ij})$ be the matrix with entries $p_{ij}$. The $\alpha$ is of the form: $(w_1,\dots,w_n)P(w'_1,\dots,w'_n)^t$.

Let $E=\k\op E_1\op E_2$ be a graded Frobenius algebra defined by an invertible matrix $M$ as above. We can view $E$ as a 2-homogeneous algebra. In fact, let $\mu:E\ot E\to E$ be the multiplication of $E$, and let $R=\{\alpha\in E_1\ot E_1|\mu(\alpha)=0\}$. Then we have $E\cong T(E_1)/(R)$.
Let $V=E_1^*$ be the dual space of $E_1$. Let $f=(x_1^*,\dots,x_n^*)M(x_1^*,\dots,x_n^*)^t\in V\ot V$. For any $x,y\in E_1$, $f(x\ot y)=\mathbf{a}M\mathbf{b}^t=\langle x,y\rangle$, where $\mathbf{a}$ and $\mathbf{b}$ are the coordinates with respect to $x$ and $y$ on the basis $\{x_1,\dots,x_n\}$. Then we see $f\in R^\bot$. So, $E^!=T(V)/(R^\bot)=T(V)/(f)$. By \cite[Theorem 0.1]{Z}, $E^!$ is a Koszul algebra. Hence $E$ is a Koszul algebra. Since $E$ is graded Frobenius, we have $E^!$ is AS-Gorenstein of global dimension 2 \cite[Proposition 5.10]{Smi}.

Summarizing the above argument, we have the following properties (see also \cite[Theorem 3]{DV}).

\begin{prop}
Let $M$ be an invertible $n\times n$-matrix, let $E=\k\op E_1\op E_2$ be the graded Frobenius algebra defined by $M$. Then:

{\rm(i)} $E$ is a Koszul algebra;

{\rm(ii)} $A=T(V)/(f)$ is AS-Gorenstein algebra of global dimension 2, where $V=E_1^*$ and $f=(x^*_1,\dots,x^*_n)M(x^*_1,\dots,x^*_n)^t\in V\ot V$;

{\rm(iii)} any AS-Gorenstein algebra of global dimension 2 is defined by an invertible matrix.
\end{prop}

On the graded Frobenius algebra $E$ defined by an invertible matrix $M$, any linear map $\delta_1:E_1\to E_2$ defines a differential $\delta_E$ on $E$. For any element $\theta\in E_2$, $(E,\delta_E,\theta)$ is a curved differential graded algebra. By \cite[Proposition 4.1, Ch. 5]{PPo}, we see the following:

\begin{prop}
Let $A=T(V)/(f)$ where $f=(x^*_1,\dots,x^*_n)M(x^*_1,\dots,x^*_n)^t\in V\ot V$ and $M$ is an invertible matrix. For any elements $v\in V$ and $k\in \k$, the filtered algebra $U=T(V)/(f-v-k)$ is a deformation of $A$.
\end{prop}

To simplify the notions, in the rest of this section, we interchange $E_1$ with $V=E_1^*$, and we fix a basis $\{x_1,\dots,x_n\}$ of $V$.

Since $A=T(V)/(f)$ is Koszul, the Nakayama automorphism of $A$ is determined by the Nakayama automorphism of $E$. If a matrix $C$ defines an automorphism of a finite dimensional vector space $W$ then its dual automorphism on $W^*$ is defined by $C^t$. Thus by Lemma \ref{lem12}, we obtain the following:

\begin{prop}\label{mprop1}
Let $A=T(V)/(f)$ where $f=(x_1,\dots,x_n)M(x_1,\dots,x_n)^t\in V\ot V$ and $M$ is an invertible matrix. Then we have $\Ext^i_{A^e}(A,A\ot A)=0$ for $i\neq 2$, and $$\Ext_{A^e}^2(A,A^e)\cong {}_1A_\zeta(-2),$$ where the Nakayama automorphism $\zeta$ is defined by $\zeta(y)=-(x_1,\dots,x_n)M^tM^{-1}\mathbf{k}^t$ and $y=(x_1,\dots,x_n)\mathbf{k}^t\in V$.
\end{prop}

As a consequence, we we obtain results of Berger (see \cite[Propositions 3.4 and 6.3]{Ber}):

\begin{cor}
Let $A=T(V)/(f)$ where $f=(x_1,\dots,x_n)Q(x_1,\dots,x_n)^t\in V\ot V$ and $Q$ is an $n\times n$ matrix. Then $A$ is Calabi-Yau of dimension 2 if and only if $Q$ is invertible and anti-symmetric.
\end{cor}
\proof The fact that $A$ is Calabi-Yau implies that $A$ is AS-Gorenstein. Hence $Q$ must be invertible. Then the result follows from Proposition \ref{mprop1}. \qed

We may write down the Nakayama automorphism of a deformation of $A$ even more precisely than it in Section 2. Let $U=T(V)/(f-v-k)$ be a deformation of $A=T(V)/(f)$ with $f=(x_1,\dots,x_n)Q(x_1,\dots,x_n)^t\in V\ot V$. Assume $v=(x_1,\dots,x_n)\mathbf{s}^t$ with $\mathbf{s}=(s_1,\dots,s_n)$. The projective resolution of the bimodule ${}_UU_U$ formed in Section 2 reads as follows: (also cf. \cite{Ber}):
$$\xymatrix{
  0 \ar[r]&U\ot U\ar[rr]^{\cdot^r\mathbf{x}Q+\cdot^l\mathbf{x}Q^t-\cdot\mathbf{s}}&&(U\ot U)^{1\times n}\ar[rr]^{\cdot^r\mathbf{x}-\cdot^l\mathbf{x}}&&U\ot U\ar[r]&U\ar[r]&0,}$$
where $(a\ot b)\cdot^r c=ac\ot b$ and $(a\ot b)\cdot^lc=a\ot cb$ for $a,b,c\in U$. After applying the functor $\Hom_{U^e}(-,U\ot U)$ to the forgoing resolution we get
\begin{equation}\label{mtag1}
    \xymatrix{0\ar[r]&U\ot U\ar[rr]^{\cdot^r\mathbf{x}-\cdot^l\mathbf{x}}&&(U\ot U)^{1\times n}\ar[rr]^{\cdot^rQ\mathbf{x}^t+\cdot^lQ^t\mathbf{x}^t-\cdot\mathbf{s}^t}&&U\ot U.}
\end{equation}
By Lemma \ref{lem21}, we see that the above complex is exact except at the final position, and the cohomology of the final position is ${}_1U_{\xi}$. Hence there is a surjective bimodule morphism $\pi:U\ot U\to {}_1U_{\xi}$ such that $\pi\circ(\cdot^rQ\mathbf{x}^t+\cdot^lQ^t\mathbf{x}^t-\cdot\mathbf{a}^t)=0$.
Similar computations to those in Section 2 yield:
\begin{equation}\label{mtag2}
    \sum_{i,j}a_iq_{ij}x_j\xi(b_i)+\sum_{i,j}a_i\xi(x_jq_{ji}b_i)-\sum_{i}s_ia_i\xi(b_i)=0,
\end{equation} for any $a_1,\dots,a_n, b_1\dots,b_n\in U$.
By Lemma \ref{lem21}, we have $\xi(x_i)=\zeta(x_i)+\lambda_i$ for $i=1,\dots,n$ with $\lambda_i\in \k$. Let $\lambda=(\lambda_1,\dots,\lambda_n)$. In the equation (\ref{mtag2}), let $b_1=\cdots=b_n=1$. Then we obtain:
\begin{equation}
   \nonumber \mathbf{a}Q\mathbf{x}^t-\mathbf{a}Q^t(Q^{-1})^tQ\mathbf{x}^t+\mathbf{a}Q^t\lambda^t-\mathbf{a}\mathbf{s}^t=0,
\end{equation} where $\mathbf{a}=(a_1,\dots,a_n)$ and $Q=(q_{ij})$.
Equivalently, we have:
\begin{equation}\label{mtag3}
   \mathbf{a}Q^t\lambda^t-\mathbf{a}\mathbf{s}^t=0.
\end{equation}
Since $\mathbf{a}$ is arbitrary, the equation (\ref{mtag3}) yields: $$Q^t\lambda^t=\mathbf{s}^t.$$ Hence $$\lambda=\mathbf{s}Q^{-1}.$$
Summarizing the above argument we obtain the following theorem.

\begin{thm} \label{mthm1} Let $Q$ be an invertible $n\times n$ matrix, $A=T(V)/(f)$ with $f=\mathbf{x}Q\mathbf{x}^t$ and $U=T(V)/(f-v-k)$ a deformation of $A$, where $\mathbf{x}=(x_1,\dots,x_n)$. Assume $v=\mathbf{x}\mathbf{s}^t$ with $\mathbf{s}=(s_1,\dots,s_n)$. Then $\Ext^i_{U^e}(U,U\ot U)=0$ for $i\neq 2$, and $$\Ext^2_{U^e}(U,U\ot U)\cong {}_1U_{\xi},$$ where the automorphism $\xi$ is defined by $\xi(\mathbf{x})=-\mathbf{x}Q^tQ^{-1}+\mathbf{s}Q^{-1}$.
\end{thm}

\begin{rem}
{\rm Note that the presentation of the Nakayama isomorphism in Theorem \ref{mthm1} is different from the one in Proposition \ref{mprop2}. The reason is that we chose different bases in these two situations.}
\end{rem}

As a special case, we have the following corollary (cf. \cite{Ber}).

\begin{cor}
Let $U$ be as above. Then $U$ is Calabi-Yau if and only if $Q$ is an anti-symmetric matrix and $v=0$.
\end{cor}
\proof The sufficiency follows from Theorem \ref{mthm1}. If $U$ is Calabi-Yau, then the Nakayama automorphism $\xi$ in Theorem \ref{mthm1} is inner. However, since $A$ is a domain (cf. \cite{DV}), $U$ has no nontrivial inner automorphism. Hence $\xi=id$, and so $Q$ is anti-symmetric and $v=0$.\qed

\section{Deformations of AS-Gorenstein algebras of global dimension 3}

It is well known that an AS-Gorenstein algebra of global dimension 3 is $N$-Koszul $(N\ge2)$. Throughout this section, we assume that $A=T(V)/(R)$ is an AS-Gorenstein $N$-Koszul ($N\ge2$) algebra of global dimension 3. As before, we assume that $\dim V=n$, and fix a basis $\{x_1,\dots,x_n\}$. Let $\{x^*_1,\dots,x_n^*\}$ be the dual basis of $V^*$. Since $A$ is AS-Gorenstein, $\dim R=n$. We fix a basis of $R$, say, $\{r_1,\dots,r_n\}$.

Let $C_0=\k, C_{-1}=V, C_{-2}=R, C_{-3}=R\ot V\cap V\ot R$ and $C=C_{-3}\op C_{-2}\op C_{-1}\op \k$. As $A$ is AS-Gorenstein, we also have $\dim C_{-3}=1$. We fix a basis $z$ of $C_{-3}$. As originally suggested in \cite{AS}, we may express $z\in R\ot V$ in the form $$z=\mathbf{r}Q^{(1)}\mathbf{x}^t,$$ where $Q^{(1)}$ is an $n\times n$ matrix and $\mathbf{r}=(r_1,\dots,r_n)$. On the other hand, $z$ is an element in $V\ot R$, so there is an $n\times n$ matrix $Q^{(2)}$ such that $$z=\mathbf{x}Q^{(2)}\mathbf{r}^t.$$ Define a coproduct on $C$ as follows:
$\Delta(v)=1\ot v+v\ot 1$ for $v\in V$, $\Delta(r)=1\ot r+r\ot 1$ for $r\in R$, and $\Delta(z)=1\ot z+z\ot1+\sum_{i=1}^nq^{(1)}_{ij}r_i\ot x_j+\sum_{i=1}^nq^{(2)}_{ij}x_i\ot r_j$. We may define a counit $\varepsilon$ on $C$ in an obvious way so that $(C,\Delta,\varepsilon)$ is a graded coalgebra. Note that $C$ in general is not a subcoalgebra of the tensor coalgebra of $V$. Let us consider the Yoneda Ext-algebra $E(A)=\op_{i\ge0}\Ext_A^i(\k_A,\k_A)$ of $A$. By \cite[Proposition 3.1]{BM} or \cite[Theorem 9.1]{GMMZ}, we see that $E(A)$ is the dual algebra of the graded coalgebra $C$.

Since $A$ is AS-Gorenstein, $E(A)$ is a graded Frobenius algebra. Now let us check the multiplication of $E(A)$. Let $x_1^*,\dots,x_n^*$ be the basis of $V^*$ dual to $x_1,\dots,x_n$, and $r_1^*,\dots,r_n^*$ be the basis of $R^*$ dual to $r_1,\dots,r_n$. For any element $\alpha,\beta\in V^*\cong \Ext_A^1(\k_A,\k_A)$, assume $\alpha=a_1x_1^*+\cdots+a_nx_n^*$ for some $\mathbf{a}=(a_1,\dots,a_n)\in \k^n$, and $\beta=b_1x_1^*+\cdots+b_nx_n^*$ for some $\mathbf{b}=(b_1,\dots,b_n)\in \k^n$. For any $\gamma \in R^*\cong \Ext_A^2(\k_A,\k_A)$, we assume $\gamma=c_1r^*_1+\cdots+c_nr^*_n$ for some $\mathbf{c}=(c_1,\dots,c_n)\in \k^n$. As $E(A)$ is the dual algebra of the coalgebra $C$, we have $$\begin{array}{ccl}
  (\gamma\cdot\alpha)(z)&=&(\gamma\ot \alpha)\Delta(z)\\
  &=&\displaystyle(\gamma\ot \alpha)(1\ot z+z\ot 1+\sum_{i,j=1}^nq^{(1)}_{ij}r_i\ot x_j+\sum_{i,j,k=1}^nq^{(1)}_{ij}p^{(i)}_{sk}x_s\ot (x_k\ot x_j)) \\
  &=&\displaystyle \sum_{i,j=1}^nc_iq^{(1)}_{ij}a_j\\
  &=& \mathbf{c}Q^{(1)}\mathbf{a}^t.
\end{array}$$ Hence we obtain: $\gamma\cdot\alpha=\mathbf{c}Q^{(1)}\mathbf{a}^tz^*$. Similarly, we have: $\alpha\cdot\gamma=\mathbf{a}Q^{(2)}\mathbf{c}^tz^*$.

If $N=2$, then $R\in V\ot V$. So, there are $n\times n$ matrices $P^{(1)},\dots,P^{(n)}$ such that $$r_1=\mathbf{x}P^{(1)}\mathbf{x}^t,\dots,r_n=\mathbf{x}P^{(n)}\mathbf{x}^t$$ respectively, where $\mathbf{x}=(x_1,\dots,x_n)$ comes from the basis. Now, we have $\alpha\cdot\beta\in R^*\cong\Ext^2_A(\k_A,\k_A)$. Then for all $i=1,\dots,n$, we have
$$\begin{array}{ccl}
  (\alpha\cdot\beta)(r_i)&=&(\alpha\ot \beta)\Delta(r_i)\\
  &=&\displaystyle(\alpha\ot \beta)(1\ot r_i+r_i\ot 1+\sum_{j,k=1}^np^{(i)}_{jk}x_j\ot x_k) \\
  &=&\displaystyle \sum_{j,k=1}^na_jp^{(i)}_{jk}b_k\\
  &=& \mathbf{a}P^{(i)}\mathbf{b}^t.
\end{array}
$$ Hence $\alpha\cdot\beta=\sum_{i=1}^n\mathbf{a}P^{(i)}\mathbf{b}^tr_i^*$.

If $N\ge3$, then $\alpha\cdot\beta=0$ (see \cite{BM,GMMZ}).

Since $E(A)$ is graded Frobenius, we have that for any $\alpha\in E^1(A)$ there is an element $\gamma\in E^2(A)$ such that $\gamma\cdot\alpha\neq 0$, which implies that $Q^{(1)}$ is invertible. Similarly, we have that $Q^{(2)}$ is invertible.

In summary, we have the following properties.

\begin{prop}\label{mprop4}
With the notations as above,

{\rm(i)} the matrices $Q^{(1)}$ and $Q^{(2)}$ are invertible;

{\rm(ii)} for $\alpha=a_1x_1^*+\cdots+a_nx_n^*,\beta=b_1x_1^*+\cdots+b_nx_n^*\in V^*\cong \Ext_A^1(\k_A,\k_A)$ and $\gamma=c_1r^*_1+\cdots+c_nr^*_n\in R^*\cong \Ext_A^2(\k_A,\k_A)$, we have
$$\gamma\cdot\alpha=\mathbf{c}Q^{(1)}\mathbf{a}^tz^*,\text{\ }\alpha\cdot\gamma=\mathbf{a}Q^{(2)}\mathbf{c}^tz^*,$$
where $\mathbf{a}=(a_1,\dots,a_n)\in \k^n$ and $\mathbf{c}=(c_1,\dots,c_n)\in \k^n$.

If $N=2$, $\displaystyle\alpha\cdot\beta=\sum_{i=1}^n\mathbf{a}P^{(i)}\mathbf{b}^tr_i^*$; if $N\ge3$, $\alpha\cdot\beta=0$.
\end{prop}

Through the isomorphism $E(A)\cong C^*$, we may view the elements in $E(A)$ as linear maps from $C$ to $\k$. Define a bilinear form $\langle\ ,\ \rangle:E(A)\times E(A)\to \k$ by $\langle\vartheta,\varsigma\rangle=(\vartheta\cdot\varsigma)(z)$ for all $\vartheta,\varsigma\in E(A)$. It is easy to see that $\langle\ ,\ \rangle$ is nondegenerate. Let $\alpha$ and $\gamma$ be as in Proposition \ref{mprop4}. We have $\langle\gamma,\alpha\rangle=\mathbf{c}Q^{(1)}\mathbf{a}^t$ and $\langle\alpha,\gamma\rangle=\mathbf{a}Q^{(2)}\mathbf{c}^t$. Since $\mathbf{a}Q^{(2)}\mathbf{c}^t=\mathbf{c}{Q^{(2)}}^t\mathbf{a}^t=\mathbf{c}Q^{(1)}{Q^{(1)}}^{-1}{Q^{(2)}}^t\mathbf{a}^t$, we see that the Nakayama automorphism $\phi$ of $E(A)$ acting on $\Ext^1_A(\k,\k)\cong V^*$ yields $$\phi(x_1^*,\dots,x_n^*)=(x_1^*,\dots,x^*_n){Q^{(1)}}^{-1}{Q^{(2)}}^t.$$

Now applying Lemma \ref{lem12}, we have the following result.
\begin{prop}\label{mprop5}
We have $\Ext^3_{A^e}(A,A\ot A)\cong {}_1A_{\zeta}(N+1)$, where $\zeta$ acts on $A_1=V$ by
$$  \zeta(x_1,\dots,x_1)=(x_1,\dots,x_n){Q^{(1)}}^t{Q^{(2)}}^{-1}.$$
\end{prop}

If $A$ is Calabi-Yau of dimension 3, then $E(A)$ is actually a symmetric algebra (\cite[Proposition A.5.2]{Bo}). Hence Proposition \ref{mprop4} implies that $Q^{(1)}={Q^{(2)}}^t$. Conversely, if $Q^{(1)}={Q^{(2)}}^t$, then Proposition \ref{mprop5} says that $A$ is Calabi-Yau. Therefore we have

\begin{cor} \label{mcor1} $A$ is Calabi-Yau if and only if $Q^{(1)}={Q^{(2)}}^t$.
\end{cor}

Let $w\in V^{\ot m}$ and $\alpha\in V^*$, we write $[\alpha w]=(\alpha\ot 1\ot\cdots\ot1)(w)$ and $[w\alpha]= (1\ot\cdots\ot1\ot\alpha)(w)$. We may extend this notion to nonhomogeneous case, that is, if $w=w_1+w_2+\cdots+w_m$ with $w_i\in V^{\ot i}$ for all $i=1,\dots,m$, then $[\alpha w]=[\alpha w_1]+\cdots+[\alpha w_m]$. Similarly, we have $[w\alpha]$. Recall that an element $w\in V^{\ot m}$ is called a {\it superpotential} of degree $m$ if $[\alpha w]=[w\alpha]$ for all $\alpha\in V^*$ \cite{Gin,Bo,VdB2}. In general, an element $w\in T^{\ge1}(V)$ is called a {\it potential} if $[\alpha w]=[w\alpha]$ for all $\alpha\in V^*$ \cite{BT}. Given an element $\alpha\in V^*$, the {\it partial derivation} of a (super)potential $w$ by $\alpha$ is defined to be $\partial_{\alpha}(w)=[\alpha w]$.

\begin{prop}\label{mprop6} If $A=T(V)/(R)$ is Calabi-Yau, then any nonzero element in $R\ot V\cap V\ot R$ is a superpotential.
\end{prop}
\proof It suffices to prove that the element $z$ that we chose at the beginning of this section as a fixed basis of $R\ot V\cap V\ot R$ is a superpotential. Recall that $z=\sum_{ij}q^{(1)}_{ij}r_i\ot x_j=\sum_{ij}q^{(2)}_{ij}x_i\ot r_j$. For all $k=1,\dots,n$, we have $[x_k^*z]=\sum_{j}q^{(2)}_{kj}r_j$ and $[zx_k^*]=\sum_{i}q^{(1)}_{ik}r_i$. By Corollary \ref{mcor1}, we have $(q^{(1)}_{1k},\dots,q^{(1)}_{nk})=(q^{(2)}_{k1},\dots,q^{(2)}_{kn})$. Hence $[x_k^*z]=[zx_k^*]$ for all $k=1,\dots,n$. Therefore $z$ is a superpotential. \qed

The following result was first proved by Bocklandt in \cite{Bo} under the more general assumption that $A$ is obtained from a finite quiver. Bocklandt proved the result by a sophisticated analysis of the minimal graded projective resolution of the trivial $A$-module $A_0$. However, in the connected case, we have a quite simpler proof.

\begin{thm}{\rm(cf. \cite{Bo})} \label{mthm4} Let $A=T(V)/(R)$ be a Calabi-Yau algebra of dimension 3. Then $A$ is defined by a superpotential.

More precisely, for any nonzero element $w\in R\ot V\cap V\ot R$, we have $A=T(V)/(\partial_\alpha(w):\alpha\in V^*)$.
\end{thm}
\proof Note that $A$ is $N$-Koszul. It suffices to prove $A=T(V)/(\partial_\alpha(z):\alpha\in V)$ where $z$ is the fixed basis of $R\ot V\cap V\ot R$ that we selected at the beginning of this section. For $k=1,\dots,n$, let $r'_k=\partial_{x_k^*}(z)=[zx_k^*]=\sum_{i}q^{(1)}_{ik}r_i$. Then $$(r'_1,\dots,r'_n)=(r_1,\dots,r_n)Q^{(1)}.$$ Since $Q^{(1)}$ is invertible, $r'_1,\dots,r'_n$ is a basis of $R$. Hence the result follows. \qed

The deformations of general $N$-Koszul algebras are much more complicated than those of Koszul algebras (see \cite{FV,BGi}). We have some further efforts to make in order to consider the Calabi-Yau property of a deformation of an AS-Gorenstein algebra of dimension 3.

Let $U=T(V)/(r+\alpha_1(r)+\cdots+\alpha_{N}(r):r\in R)$ be a deformation of $A$. To determine when $U$ is Calabi-Yau, we have to construct a projective resolution for the bimodule ${}_UU_U$. We need some additional notations in the following discussions. Let $\mu:U\ot U\to U$ be the multiplication map of $U$. Denote by $\mu^0=id$, $\mu^1=\mu$, and $\mu^k=\mu(\mu^{k-1}\ot1)$ for $k\ge2$. Let $\tau:V\to U$ be the natural inclusion map. We denote $\tau^k=\tau^{\ot k}$ for $k\ge2$.
Let us consider the sequence:
{\small \begin{equation}\label{mtag11}
    0\longrightarrow U\ot(V\ot R\cap R\ot V)\ot U\overset{D^{-3}}\longrightarrow U\ot R\ot U\overset{D^{-2}}\longrightarrow U\ot V\ot U\overset{D^{-1}}\longrightarrow U\ot U\overset{\mu}\longrightarrow U\longrightarrow0,
\end{equation}}
where the morphisms are given as following:
$$\begin{array}{ccl}
    D^{-3}&=& \mu(1\ot\tau)\ot 1^{\ot N}\ot 1-1\ot 1^{\ot N}\ot\mu(\tau\ot1) \\
    &&+1\ot(\alpha_1\ot 1)\ot1-1\ot(1\ot \alpha_1)\ot1,\\[3mm]
    D^{-2}&=& \displaystyle\sum^{i+j=N-1}_{i,j\ge0}\mu^i(1\ot\tau^i)\ot1\ot\mu^j(\tau^j\ot1)\\
    &&\displaystyle+\sum_{k=1}^{N-1}\sum_{i,j\ge0}^{i+j=N-k-1}(\mu^i(1\ot\tau^i)\ot1\ot\mu^j(\tau^j\ot1))(1\ot\alpha_k\ot1),\\[3mm]
    D^{-1}&=& \mu(1\ot\tau)\ot1-1\ot\mu(\tau\ot1).
  \end{array}
$$
\begin{lem}\label{lem1} The sequence (\ref{mtag11}) is exact, hence it is a projective resolution of ${}_UU_U$.
\end{lem}
\proof We first show that the sequence (\ref{mtag11}) is a complex. Note that the morphisms in the sequence are clearly $U$-bimodule morphisms. In the following computation, we make use of the identities in \cite[Theorem 1.1]{FV}.
$$\begin{array}{cl}
    D^{-2}D^{-3}&=\mu^N(1\ot\tau^N)\ot1\ot1-1\ot1\ot\mu^N(\tau^N\ot1)\\
    &\displaystyle+\sum_{i,j\ge0}^{i+j=N-1}(\mu^i(1\ot\tau^i)\ot1\ot\mu^j(\tau^j\ot1))(1\ot\alpha_1\ot1\ot1) \\
    &\displaystyle-\sum_{i,j\ge0}^{i+j=N-1}(\mu^i(1\ot\tau^i)\ot1\ot\mu^j(\tau^j\ot1))(1\ot1\ot\alpha_1\ot1)\\
    &\displaystyle+\sum_{k=1}^{N-1}\sum_{i,j\ge0}^{i+j=N-k-1}(\mu^{i+1}(1\ot\tau^{i+1})\ot1\ot\mu^j(\tau^j\ot1))(1\ot1\ot\alpha_k\ot1)\\
    &\displaystyle-\sum_{k=1}^k\sum_{i,j\ge0}^{i+j=N-k-1}(\mu^i(1\ot\tau^i)\ot1\ot\mu^{j+1}(\tau^{j+1}\ot1))(1\ot\alpha_k\ot1\ot1)\\
    &\displaystyle-\sum_{k=1}^{N-1}\sum_{i,j\ge0}^{i+j=N-k-1}(\mu^i(1\ot\tau^i)\ot1\ot\mu^{j}(\tau^{j}\ot1))
    (1\ot\alpha_k(1\ot\alpha_1-\alpha_1\ot1)\ot1).
  \end{array}
$$
It is sufficient to show that $D^{-2}D^{-3}(1\ot z\ot1)=0$. Recall that $z=\sum_{i=1}^nq^{(1)}_{ij}r_i\ot x_j=\sum_{i=1}^nq^{(2)}_{ij}x_i\ot r_j$.
Then: $$(\mu^{N-1}\tau^N\ot1)(z)=-((\mu^{N-2}\tau^{N-1}\alpha_1+\mu^{N-3}\tau^{N-2}\alpha_{2}+\cdots+\tau\alpha_{N-1}+\alpha_N)\ot1)(z),$$ and:
$$(1\ot \mu^{N-1}\tau^N)(z)=-(1\ot(\mu^{N-2}\tau^{N-1}\alpha_1+\mu^{N-3}\tau^{N-2}\alpha_{2}+\cdots+\tau\alpha_{N-1}+\alpha_N))(z).$$
Keeping in mind that the morphisms in the above calculation would act on the element $1\ot z\ot1$, we have:
$$\mu^N(1\ot\tau^N)\ot1\ot1=-(\mu^{N-1}(1\ot\tau^{N-1}\alpha_1)+\mu^{N-2}(1\ot\alpha_2)+\cdots+\mu(1\ot\alpha_N))\ot1\ot1,$$ and:
$$-1\ot1\ot\mu^N(\tau^N\ot1)=1\ot1\ot(\mu^{N-1}(\tau^{N-1}\alpha_1\ot1)+\mu^{N-2}(\alpha_2\ot1)+\cdots+\mu(\alpha_N\ot1)).$$
Moreover, $$\begin{array}{cl}
             &\displaystyle\sum_{k=1}^{N-1}\sum_{i,j\ge0}^{i+j=N-k-1}(\mu^i(1\ot\tau^i)\ot1\ot\mu^{j}(\tau^{j}\ot1))
    (1\ot\alpha_k(1\ot\alpha_1-\alpha_1\ot1)\ot1)\\
    &\displaystyle=\sum_{k=1}^{N-1}\sum_{i,j\ge0}^{i+j=N-k-1}(\mu^i(1\ot\tau^i)\ot1\ot\mu^{j}(\tau^{j}\ot1))
    (1\ot(1\ot\alpha_{k+1}-\alpha_{k+1}\ot1)\ot1)\\
    &\displaystyle=\sum_{k=2}^{N}\sum_{i\ge0,j\ge1}^{i+j=N-k}(\mu^i(1\ot\tau^i)\ot1\ot\mu^{j}(\tau^{j}\ot1))
    (1\ot(1\ot\alpha_{k}-\alpha_{k}\ot1)\ot1).
           \end{array}
$$
Comparing the components in the three equations above with the components of $D^{-3}D^{-2}$, we see that $D^{-3}D^{-2}=0$.

Similarly, we can check that $D^{-2}D^{-1}=0$. The equation $\mu D^{-1}=0$ holds obviously. Hence (\ref{mtag11}) is a complex. We may define a filtration on the complex (\ref{mtag11}) as we did in Section \ref{gendef} so that it is a complex of filtered $U$-bimodules: the filtration on $U\ot U$ is just the filtration induced by that of $U$; $F_{i}(U\ot V\ot U)=0$ for $i<1$, and $\displaystyle F_{1+k}(U\ot V\ot U)=\sum_{i+j=k}F_iU\ot V\ot F_jU$ for $k\ge0$; $F_{i}(U\ot R\ot U)=0$ for $i<N$, and $\displaystyle F_{N+k}(U\ot R\ot U)=\sum_{i+j=k}F_iU\ot R\ot F_jU$; $\displaystyle F_{i}(U\ot (R\ot V\cap V\ot R)\ot U)=0$ for $i<N+1$, and $\displaystyle F_{N+1+k}(U\ot (R\ot V\cap V\ot R)\ot U)=\sum_{i+j=k}F_iU\ot (R\ot V\cap V\ot R)\ot F_jU$. The first level of the spectral sequence of the complex (\ref{mtag11}) with the filtration given above is just the bimodule Koszul complex of $A$ as constructed in \cite[Theorem 4.4]{BM}. Hence the sequence (\ref{mtag11}) is exact. \qed

Similar to Lemma \ref{lem21}, we also have the following result obtained from the projective resolution (\ref{mtag11}) of $U$. The proof is very similar to that of Lemma \ref{lem21}, so we omit it.

\begin{lem} \label{lem31} Let $A$ be an $N$-Koszul AS-Gorenstein algebra of global dimension 3. Let $U$ be a deformation of $A$. Assume that the Nakayama automorphism of $A$ is $\zeta$. Then $\Ext^i_{U^e}(U,U\ot U)=0$ for $i\neq 3$, and $\Ext^3_{U^e}(U,U\ot U)\cong {}_1U_\xi$ for some filtration-preserving automorphism $\xi$ such that $gr(\xi)=\zeta$.
\end{lem}

Applying the functor $\Hom_{U^e}(-,U\ot U)$ to the projective resolution (\ref{mtag11}) of ${}_UU_U$, we get the following complex:
{\small \begin{equation}\label{mtag12}
    0\longrightarrow U\ot U\overset{\widehat{D}^{-1}}\longrightarrow U\ot V^*\ot U\overset{\widehat{D}^{-2}}\longrightarrow U\ot R^*\ot U\overset{\widehat{D}^{-3}}\longrightarrow U\ot(R\ot V\cap V\ot R)^*\ot U.
\end{equation}}
We only write down explicitly the morphism $\widehat{D}^{-3}$. Write $\psi$ for the morphism $1\ot \alpha_1-\alpha_1\ot1: R\ot V\cap V\ot R\longrightarrow R$, and let $\psi^*$ be the dual map of $\psi$. For $\vartheta\in R^*$, we have: $$\widehat{D}^{-3}(1\ot\vartheta\ot1)=\sum_{i=1}^n1\ot x^*_i\vartheta\ot \tau(x_i)-\sum_{i=1}^n\tau(x_i)\ot \vartheta x_i^*\ot 1-1\ot \psi^*(\vartheta)\ot 1,$$
where the products $x^*_i\vartheta$ and $\vartheta x_i^*$ are the Yoneda products in $E(A)$.

\begin{prop}\label{mprop13} Keep the preceding notions. Assume $\psi^*(r_i)=k_iz^*$ for $i=1,\dots,n$, and let $\mathbf{k}=(k_1,\dots,k_n)\in \k^n$. Then the Nakayama automorphism of $U$ is defined by the linear map
$\chi:V\to V\op \k$, with $\chi(x_1,\dots,x_n)=(x_1,\dots,x_n){Q^{(1)}}^t{Q^{(2)}}^{-1}+\mathbf{k}{Q^{(2)}}^{-1}$.
\end{prop}
\proof We only need to compute the cohomology at the final position of the complex (\ref{mtag12}). Since the computation is very similar to the one in Proposition \ref{mprop2}, we omit it. \qed

From now on, we will further assume that $A$ is Calabi-Yau.

\begin{cor}\label{mcor2} Let $U=T(V)/(r+\alpha_1(r)+\cdots+\alpha_{N-1}(r):r\in R)$ be an augmented deformation of $A$.
Then $U$ is Calabi-Yau if and only if $E(U)=\Ext_U^*(\k,\k)$ is a symmetric algebra of length 3.
\end{cor}
\proof It suffices to prove the ``if'' part. \cite[Theorem 2.1]{FV} shows that $E(A)=\Ext_A^*(\k,\k)$ is an $A_\infty$-algebra, and \cite[Theorem 2.3]{FV} says that $E(U)$ is the cohomolgy algebra of the $A_\infty$-algebra $E(A)$. Since $\dim\Ext^3_A(\k,\k)=1$, the action of the differential $m_1$ on $\Ext^2_A(\k,\k)$ is zero. On the other hand, the restrict of $m_1$ to $\Ext^2_A(\k,\k)=R^*$ is exactly the map $\psi^*$ as formed in the statement above. Hence Proposition \ref{mprop13} implies that the Nakayama automorphism is the identity map. Therefore $U$ is Calabi-Yau. \qed

Let $A=T(V)/(R)$ be an $N$-Koszul Calabi-Yau algebra of dimension 3. Assume that $U=T(V)/(r+\alpha_1(r)+\cdots+\alpha_{N}(r):r\in R)$ is a deformation of $A$, and that $U$ is Calabi-Yau. Assume further that $(\alpha_1\ot 1 -1\ot\alpha_1)(z)=0$. Then by \cite[Theorem 1.1]{FV} or \cite[Proposition 3.6]{BGi}, we have $(\alpha_k\ot 1 -1\ot\alpha_k)(z)=0$ for all $k=2,\dots,N$. Hence we have
\begin{equation}\label{mtag13}
\sum_{i,j=1}^nq^{(1)}_{ij}\alpha_k(r_i)\ot x_j=\sum_{i,j=1}^nq^{(2)}_{ij}x_i\ot\alpha_k(r_j),
\end{equation}
for all $k=1,\dots,N$.

\begin{lem} \label{lem41}
If $(\alpha_1\ot 1 -1\ot\alpha_1)(z)=0$, then $w=z+(\alpha_1\ot1)(z)+\cdots+(\alpha_{N}\ot1)(z)$ is a potential.
\end{lem}
\proof It suffices to show that $g_k:=(\alpha_k\ot 1)(z)$ is a superpotential for all $k=1,\dots,N-1$ since we already know that $z$ is a superpotential. For $s=1,\dots,n$, we have $[g_k x_s^*]=\sum_{i=1}^nq^{(1)}_{is}\alpha_k(r_i)$, and $[x_s^*g_k]=\sum_{j=1}^nq^{(2)}_{sj}\alpha_k(r_j)$ by the equation (\ref{mtag13}). Since $Q^{(1)}={Q^{(2)}}^t$, we have $[x_s^*g_k]=[g_k x_s^*]$ for all $s=1,\dots,n$. That is, $g_k$ is a superpotential. \qed

The following result can be viewed as a converse of Berger-Taillefer's result \cite[Theorem 3.6]{BT}.

\begin{thm}\label{mthm5}
Let $A=T(V)/(R)$ be an $N$-Koszul Calabi-Yau algebra of dimension 3. Assume that $U$ is a deformation of $A$ which is also Calabi-Yau. If one of the following conditions is satisfied, then $U$ is defined by a potential $w\in T^{\ge1}(V)$.

{\rm(i)} $U$ is an augmented deformation of $A$;

{\rm(ii)} $A$ is a domain.
\end{thm}
\proof Assume $U=T(V)/(r+\alpha_1(r)+\cdots+\alpha_{N}(r):r\in R)$ for some linear maps $\alpha_k:R\to V^{\ot N-k}$. (i) if $U$ is augmented, then Corollary \ref{mcor2} implies that $(\alpha_1\ot 1 -1\ot\alpha_1)(z)=0$. (ii) if $A$ is a domain, then we see that there is a unique inner isomorphism of $U$. Hence the Nakayama automorphism $\xi$ of $U$ is the identity map. Then, by Proposition \ref{mprop13}, we also have $(\alpha_1\ot 1 -1\ot\alpha_1)(z)=0$. Now let $w$ be the potential given in Lemma \ref{lem41}. For $k=1,\dots,n$, $$\partial_{x^*_s}(w)=[wx_s^*]=\sum_{i=1}^nq^{(1)}_{is}r_i+\sum_{i=1}^nq^{(1)}_{is}\alpha_1(r_i)+\cdots+\sum_{i=1}^nq^{(1)}_{is}\alpha_N(r_i).$$ Let $\psi_k=r_i+\alpha_1(r_i)+\cdots+\alpha_N(r_i)$. Then we have $(\partial_{x^*_1}(w),\dots,\partial_{x^*_n}(w))=(\psi_1,\dots,\psi_n)Q^{(1)}$. Since $Q^{(1)}$ is invertible, $span\{\partial_{x^*_1}(w),\dots,\partial_{x^*_n}(w)\}=span\{r+\alpha_1(r)+\cdots+\alpha_N(r):r\in R\}$. \qed

\subsection*{Acknowledgement} The authors thank the anonymous referee for his/her valuable suggestions and comments. The work is
supported by an FWO-grant and grants from NSFC (No. 11171067), ZJNSF (No. LY12A01013), Science and Technology Department of Zhejiang Province (No. 2011R10051), and SRF for ROCS, SEM.

\bibliography{}

\end{document}